\documentclass{article}
\title{An algorithm of computing inhomogeneous differential
equations for definite integrals}
\author{Hiromasa Nakayama, Kenta Nishiyama
 \footnote{Department of Mathematics, Kobe University and JST CREST.}
}
\date{July 14, 2010}
\usepackage{amsmath,amscd,amssymb}
\usepackage{theorem}

\newtheorem{theorem}{Theorem}

\newtheorem{remark}{Remark}
\newtheorem{example}{Example}

\newtheorem{algorithm}{Algorithm}

\def\p{\partial}
\def\pd#1{ \partial_{#1} }
\def\ord{{\mathrm ord}}

\begin{document}
\maketitle

\begin{abstract}
We give an algorithm to compute inhomogeneous differential equations 
for definite integrals with parameters.
The algorithm is based on the integration algorithm for $D$-modules by Oaku.
Main tool in the algorithm is the Gr\"obner basis method in the ring of 
differential operators.
\end{abstract}

\section{Introduction}
Let us denote by 
$D = K \langle x_1, \ldots, x_n, \pd{1}, \ldots, \pd{n} \rangle$
the Weyl algebra in $n$ variables, where $K$ is $\mathbb{Q}$ or $\mathbb{C}$ and $\pd{i}$ is 
the differential operator standing for $x_i$.
We denote by 
$D' = K \langle x_{m+1}, \ldots, x_n, \pd{m+1}, \ldots, \pd{n} \rangle$
the Weyl algebra in $n-m$ variables, 
where $m \leq n$ and $D'$ is a subring of $D$.

Let $I$ be a holonomic left $D$-ideal (\cite{SST}).
The integration ideal of $I$ with respect to $x_1, \ldots, x_m$ is defined by  
the left $D'$-ideal 
$$(I + \pd{1} D + \cdots + \pd{m} D) \cap D'.$$
Oaku (\cite{OAKU}) gave an algorithm computing the integration ideal. 
This algorithm is called {\it the integration algorithm for $D$-modules}.
The Gr\"{o}bner basis method in $D$ is used in this algorithm. 

We give a new algorithm computing  
not only generators of the integration ideal $J$ but also 
$P_0 \in I$ and $P_1, \ldots, P_m \in D$ such as 
$$P = P_0 + \pd{1} P_1 + \cdots + \pd{m} P_m$$
for any generator $P \in J$.
Our algorithm is based on Oaku's one. 
We call these $P_1, \ldots, P_m$ {\it inhomogeneous parts of $P$}. 
As an important application of our algorithm, 
we can obtain inhomogeneous differential equations for definite integrals with parameters 
by using generators of the integration ideal and inhomogeneous parts. 

For example, 
we compute an inhomogeneous differential equation for  
the integral $A(x_2) = \int_{a}^{b} e^{-x_1-x_2 x_1^3} dx_1$.
This is the case of $m=1, n=2$.
The annihilating ideal of the integrand $f(x_1,x_2) = e^{-x_1-x_2 x_1^3}$ in $D$ is 
$ I = \langle \pd{1} + 1 + 3x_2x_1^2, \pd{2} + x_1^3 \rangle.$
The integration ideal of $I$ with respect to $x_1$ is  
$ J = \langle 27x_2^3 \pd{2}^2+54x_2^2\pd{2}+6x_2+1 \rangle = \langle P \rangle.$
The operator $P_1 =  - (\pd{1}^2+3\pd{1}+3)$ is an inhomogeneous part of $P$. 
We apply the operator $P$ to the integral $A(x_2)$ and obtain 
\begin{align*}
   P \cdot A(x_2) 
   &= \int_{a}^{b} \pd{1} (P_1 \cdot  e^{-x_1-x_2 x_1^3}) dx_1
   = \left [ P_1 \cdot  e^{-x_1-x_2 x_1^3} \right ]_{x_1=a}^{x_1=b} \\ 
   &= -\left[ (9x_2^2x_1^4-3x_2x_1^2-6x_2x_1+1) e^{-x_1-x_2 x_1^3} \right]_{x_1=a}^{x_1=b}.
\end{align*}
In this way, we get an inhomogeneous differential equation for  the integral $A(x_2)$. 

We will give an algorithm to compute 
inhomogeneous parts of the integration ideal and give some examples.
Other algorithms to compute differential equations for 
definite integrals are the Almkvist-Zeilberger algorithm (\cite{AZ},
\cite{Tefera}, \cite{APZ}),
the Chyzak algorithm (\cite{C2}) and the Oaku-Shiraki-Takayama algorithm (\cite{OST}).
A comparison with these algorithms are also given.

We implement our algorithms on the computer algebra system Risa/Asir (\cite{Asir}).
They are in the program package {\tt nk\_restriction.rr} (\cite{IntPackage}).
Packages {\tt Mgfun} in Maple and {\tt HolonomicFunctions} in Mathematica offers 
an analogous functionality, and are based on the Chyzak algorithm
(\cite{Mgfun}, \cite{HolonomicFunctions}).

\section{Review of the integration algorithm for $D$-modules} \label{sec:int-ideal}

We will review the integration algorithm for $D$-modules.
We define the ring isomorphism $\mathcal{F} : D \rightarrow D$ satisfying
$$\mathcal{F}(x_i) = 
\begin{cases} 
-\pd{i} ~~~(1 \leq i \leq m) \\
x_i     ~~~(m < i \leq n)
\end{cases}
,  
\mathcal{F}(\pd{i}) = 
\begin{cases} 
x_i     ~~~(1 \leq i \leq m) \\
\pd{i}  ~~~(m < i \leq n)
\end{cases}.
$$ 
This map is called {\it the Fourier transformation in $D$}.

The integration ideal of a left holonomic $D$-ideal $I$ with respect to $x_1, \ldots, x_m$ 
is defined by the left $D'$-ideal $J = (I + \pd{1} D + \cdots + \pd{m} D) \cap D'.$

\begin{algorithm} 
[Integration algorithm for $D$-modules, \cite{OAKU}, \cite{SST}]
\label{alg:int-ideal}
~~~ 
\normalfont
\begin{itemize}
\item[Input:] Generators of a holonomic left $D$-ideal $I$ and\\
a weight vector $w = (w_1, \ldots, w_m, w_{m+1}, \ldots, w_n)$ such that 
$w_1, \ldots, w_m > 0, w_{m+1} = \cdots = w_n = 0$.
\item[Output:] Generators of the integration ideal of $I$ with respect to $x_1, \ldots, x_m$.
\end{itemize}
\begin{enumerate}
\item Compute the restriction module of the left $D$-ideal $\mathcal{F}(I)$ with respect to 
      the weight vector $w$.
      The details of the computation are as follows. 
	  \begin{enumerate}
      \item Compute the Gr\"{o}bner basis of the left $D$-ideal $\mathcal{F}(I)$ 
            with respect to the monomial order $<_{(-w,w)}$. 
            Let the Gr\"{o}bner basis be $G = \{h_1, \ldots, h_l\}$.
      \item Compute the generic $b$-function $b(s)$ of $\mathcal{F}(I)$ with respect 
            to the weight vector $(-w,w)$. 
      \item If $b(s)$ has a non-negative integer root, then we set \\
            $s_0 = (\text{the maximal non-negative integer roots})$. \\   
            Otherwise, the integration ideal is 0 and finish.
      \item $m_i = \ord_{(-w,w)}(h_i)$, \\  
            $\mathcal{B}_d = \{\pd{1}^{i_1} \cdots \pd{m}^{i_m} ~|~
            i_1 w_1 + \cdots + i_m w_m \leq d \} \quad (d \in \mathbb{N}) $ , \\
            $r = \#\{(i_1, \ldots, i_m) ~|~
            i_1 w_1 + \cdots + i_m w_m \leq s_0\} = \#\mathcal{B}_{s_0}$. 
      \item $\tilde{\mathcal{B}} = \displaystyle{\bigcup_{i=1}^l}
            \{ \tilde{h}_{i\beta} := \partial^\beta h_i 
            \mid \partial^\beta \in \mathcal{B}_{s_0-m_i} \}$,\\
            $\mathcal{B}=\{h_{i\beta}:=\tilde{h}_{i\beta}|_{x_1=\cdots=x_m=0}
            \mid \tilde{h}_{i\beta} \in \tilde{\mathcal{B}} \}$. \\
            Here, $h_{i\beta} = \sum_{\p^\alpha \in \mathcal{B}_{s_0}} g_\alpha \p^\alpha \quad (g_\alpha \in D').$
      \end{enumerate}
\item Let $(D')^r$ be the left free $D'$-module with the base $\mathcal{F}^{-1}(B_{s_0})$, i.e. 
      $(D')^r = \sum_{\p^\alpha \in \mathcal{B}_{s_0}} D' x^\alpha$.
      Regard elements in $\mathcal{F}^{-1}(\mathcal{B})$ as elements in the left $D'$-module $(D')^r$. 
      In other words, $\mathcal{F}^{-1}(h_{i\beta}) = \sum_{\p^\alpha \in \mathcal{B}_{s_0}} g_\alpha x^\alpha 
      \quad (g_\alpha \in D')$ is regarded as an element in $(D')^r$.
      Let $M$ be the left $D'$-submodule in $(D')^r$ generated by $\mathcal{F}^{-1}(\mathcal{B})$.  
\item Compute the Gr\"{o}bner basis $G$ of $M$ with respect to a POT term order such that the position
      corresponds to $x^0 = 1$ is the minimum position. 
      Output $G' = G \cap D'$. This set $G'$ generates the integration ideal of $I$. 
\end{enumerate}
\end{algorithm}

We consider the following definite integral of a holonomic function $f(x_1, \ldots, x_n)$.
$$ A(x_{m+1}, \ldots, x_n) = \int_{R} f(x_1, \ldots, x_n) dx_1 \cdots dx_m, \quad R = \prod_{i=1}^m [a_i, b_i] $$
Let $I = {\rm Ann}_{D} f := \{ P\cdot f=0  \mid P \in D \}$
be the annihilating ideal of the integrand,
and $J$ be the integration ideal of $I$.
For every $p \in J$, there exist $p_1, \ldots, p_m \in D $ such that
$$ p - \sum_{i=1}^m \pd{i} p_i \in I $$
and we have
\begin{align} 
  p \cdot A(x_{m+1}, \ldots, x_n)
   &= \int_R p \cdot f dx_1 \cdots dx_m 
   = \int_R \sum_{i=1}^m (\pd{i} p_i) \cdot f dx_1 \cdots dx_m  \notag \\
   &= \sum_{i=1}^m \int_R \pd{i} (p_i \cdot f) dx_1 \cdots dx_m .\label{eq:1}
\end{align}
Therefore, if we take an integration domain such that 
the right hand side of (\ref{eq:1}) equals to zero,
we can regard the integration ideal as 
a system of homogeneous differential equations for the integral $A(x_{m+1}, \ldots, x_n)$. 
If the right hand side is not zero, the equation (\ref{eq:1}) gives
an inhomogeneous differential equations for the function A.

\section{Computing inhomogeneous parts of the integration ideal}
In this section, we give a new algorithm of computing inhomogeneous
differential equations for definite integrals.
For the purpose, we must find an explicit form $p_i$ ($1 \leq i \leq m$)
in the equation (\ref{eq:1}) in the section \ref{sec:int-ideal}.

\begin{theorem}
Let $J \subset D'$ be the integration ideal of a holonomic left $D$-ideal $I$.
For any $p \in J$, 
there exists an algorithm to compute differential operators
$p_i \in D$ $(1 \leq i \leq m)$
such that
\begin{equation} \label{eq:inhomo-op}
p - \sum_{i=1}^m \partial_i p_i \in I .
\end{equation}
\end{theorem}
{\it Proof.}
We will present an algorithm of obtaining operators $p_i$.
By applying Algorithm \ref{alg:int-ideal},
we obtain a generating set $\{g_1,\ldots,g_t \}$ of
the integration ideal of $I$.
It is sufficient to compute inhomogeneous parts for each generator $g_j$.
From the step 3 of Algorithm \ref{alg:int-ideal},
$g_j$ can be expressed as $g_j = \sum q_{ji\beta} \mathcal{F}^{-1}(h_{i\beta})$
where $q_{ji\beta} \in D$.
Then these $q_{ji\beta} \in D$ can be computed by referring the history of
the Gr\"obner basis computation in the step 3.
Therefore, we have
\begin{align*}
I \ni \sum q_{ji\beta} \mathcal{F}^{-1}(\tilde{h}_{i\beta})
&= g_j - \left ( g_j - \sum q_{ji\beta} 
                 \mathcal{F}^{-1}(\tilde{h}_{i\beta}) \right )\\
&= g_j - \sum q_{ji\beta} \left (\mathcal{F}^{-1}(h_{i\beta}) 
                             -  \mathcal{F}^{-1}(\tilde{h}_{i\beta}) \right )\\
&= g_j - \sum q_{ji\beta} 
          \left (\mathcal{F}^{-1}(\tilde{h}_{i\beta}|_{x_1= \cdots =x_m=0})
                             -  \mathcal{F}^{-1}(\tilde{h}_{i\beta}) \right )\\
&= g_j - \sum q_{ji\beta} \mathcal{F}^{-1}
               (\tilde{h}_{i\beta}|_{x_1= \cdots =x_m=0} - \tilde{h}_{i\beta}).
\end{align*}
Since each term of
$\tilde{h}_{i\beta}|_{x_1= \cdots =x_m=0} - \tilde{h}_{i\beta}$ 
can be divided from the left by either of $x_1,\ldots,x_m$,
each term of
$\mathcal{F}^{-1}(\tilde{h}_{i\beta}|_{x_1=\cdots=x_m=0}-\tilde{h}_{i\beta})$
can be divided from the left by either of
$\partial_1,\ldots,\partial_m$.
Thus we can rewrite
$$ \sum q_{ji\beta} \mathcal{F}^{-1}(\tilde{h}_{i\beta}|_{x_1= \cdots =x_m=0}
        - \tilde{h}_{i\beta})
  = \sum_{i=1}^m \partial_i p_{ij}.$$

Let us present our algorithm.

\begin{algorithm} \normalfont \label{alg:int-inhomo}
~~~ 
\begin{itemize}
\item[Input:] Generators of a holonomic left ideal $I \subset D$ and \\
a weight vector $w = (w_1, \ldots, w_m, w_{m+1}, \ldots, w_n)$ such that 
$w_1, \ldots, w_m > 0, w_{m+1} = \cdots = w_n = 0$.
\item[Output:]  Generators $\{g_1,\ldots,g_t \}$ of the integration ideal of 
$I$ w.r.t. $x_1, \ldots, x_m$ and
operators $p_{ij} \in D$ satisfying 
$g_j - \sum_{i=1}^{m} \partial_i p_{ij} \in I$
for each generator $g_j$ $(1 \leq j \leq t)$.
\end{itemize}
\begin{enumerate}
\item Apply Algorithm \ref{alg:int-ideal}.
\item 
Compute $q_{ji\beta}$ satisfying
$g_j = \sum q_{ji\beta} \mathcal{F}^{-1}(h_{i\beta})$ 
by referring the history of the Gr\"obner basis computation
in the step 3 of Algorithm \ref{alg:int-ideal}.
\item Rewrite 
$R_j := g_j - \sum q_{ji\beta} \mathcal{F}^{-1}(\tilde{h}_{i\beta})$
to the form of 
$R_j = \displaystyle\sum_{i=1}^m \partial_i p_{ij}$.\\
Output $p_{ij}$.
\end{enumerate}
\end{algorithm}

\begin{example} \normalfont
[Incomplete Gauss's hypergeoemtric integral]

We set 
$$ F(x) = \int_p^q t^{b-1} (1-t)^{c-b-1} (1-xt)^{-a} dt.$$
We will compute a differential equation for the integral $F(x)$.
A holonomic ideal annihilating the integrand $f(x,t)=t^{b-1} (1-t)^{c-b-1} (1-xt)^{-a}$ is  
\begin{align*}
I_f = &\langle (-x^2+x)\pd{x}^2+((-t+1)\pd{t}+(-a-b-1)x+c-1)\pd{x}-ab, \\
      &(-t+1)x\pd{x}+(t^2-t)\pd{t}+(-c+2)t+b-1,(tx-1)\pd{x}+at \rangle
\end{align*}
which is obtained by using Oaku's algorithm to compute the annihilating ideal of a power of polynomials.
The generic $b$-function of $\mathcal{F}(I_f)$ with respect to the weight vector $w = (1,0)$ 
(i.e. $t$'s weight is 1 and $x$'s weight is 0) is $s (s-a+c-1)$.
We assume that $a-c+1$ is not a non-negative integer. 
Then the maximal non-negative integer root $s_0$ of $b(s)$ is 0.
Therefore, the integration ideal of $I_f$ with respect to $t$ is 
$$\langle (-x^2+x)\pd{x}^2+((-a-b-1)x+c)\pd{x}-ab \rangle = \langle P \rangle. $$
The differential equation $P \cdot g = 0$ is Gauss's hypergeometric equation. 
The inhomogeneous part of $P$ is $ \pd{t} (-t+1) \pd{x}.$
We apply $P$ to the integral $F(x)$ and obtain the inhomogeneous differential equation
$$ P \cdot \int_p^q f(x,t) dt = \int_p^q (\pd{t} (t-1) \pd{x}) \cdot f(x,t) dt = 
  \left [(t-1) \frac{\partial f}{\partial x}(x,t) \right ]_p^q.$$

We present the output for this problem by the program {\tt nk\_restriction.rr} 
on the computer algebra system Risa/Asir (\cite{Asir}). 
We use the command \\ 
{\tt nk\_restriction.integration\_ideal} to compute the integration ideal.
The option {\tt inhomo=1} make the system compute inhomogeneous parts and 
the option {\tt param = [a,b,c]} means that parameters are $a,b,c$.
The {\tt sec} shows the exhausting time of each steps.
This example and next example are executed on a Linux machine with 
Intel Xeon X5570 (2.93GHz) and 48 GB memory.

{\small
\begin{verbatim}
[1743] load("nk_restriction.rr");
[1944] I_f=[-dx^2*x^2+(-dx*a-dx*b+dx^2-dx)*x-dx*dt*t-b*a+dx*c+dx*dt-dx,
(-dx*t+dx)*x+dt*t^2+(-c-dt+2)*t+b-1,dx*t*x+a*t-dx];
[(-x^2+x)*dx^2+((-t+1)*dt+(-a-b-1)*x+c-1)*dx-b*a,
(-t+1)*x*dx+(t^2-t)*dt+(-c+2)*t+b-1,(t*x-1)*dx+a*t]
[1945] nk_restriction.integration_ideal(I_f,[t,x],[dt,dx],[1,0]|param=
[a,b,c],inhomo=1);
-- nd_weyl_gr :0.004sec(0.000623sec)
-- weyl_minipoly_by_elim :0sec(0.000947sec)
-- generic_bfct_and_gr :0.004sec(0.001922sec)
generic bfct : [[1,1],[s,1],[s-a+c-1,1]]
S0 : 0
B_{S0} length : 1
-- fctr(BF) + base :0sec(0.000277sec)
-- integration_ideal_internal :0sec(0.000499sec)
[[(-x^2+x)*dx^2+((-a-b-1)*x+c)*dx-b*a],[[[[dt,(t-1)*dx]],1]]]
\end{verbatim}
}

\end{example}

\begin{example} \normalfont
[$F(x) = \int_0^\infty e^{-t-x t^3} dt$]

We consider the integral $F(x) = \int_0^\infty e^{-t-x t^3} dt$.
A holonomic ideal annihilating the integrand $f(t,x) = e^{-t-x t^3}$ is 
$ I_f = \langle \pd{t} + 1 + 3xt^2, \pd{x} + t^3 \rangle$.
The integration ideal of $I_f$ with respect to $t$ is 
$ J = \langle 27x^3\pd{x}^2+54x^2\pd{x}+6x+1 \rangle = \langle P \rangle$.
The inhomogeneous part of $P$ is
$ -\pd{t} (\pd{t}^2+3\pd{t}+3).$
We apply $P$ to the integral $F(x)$ and obtain  
\begin{align*}
P \cdot \int_0^\infty e^{-t-x t^3} dt &= 
-\int_0^\infty (\pd{t} (\pd{t}^2+3\pd{t}+3)) \cdot  e^{-t-x t^3} dt \\
&= -\left[ (\pd{t}^2 + 3\pd{t} + 3) \cdot e^{-t-x t^3} \right]_0^\infty \\
&= -\left[ (-6xt + (1+3xt^2)^2-3-9xt^2+3) e^{-t-x t^3} \right]_0^\infty = 1 .
\end{align*}

{\small
\begin{verbatim}
[1946] load("nk_restriction.rr");
[2146] I_f=[dt+1+3*x*t^2, dx+t^3];
[dt+3*t^2*x+1,dx+t^3]
[2147] nk_restriction.integration_ideal(I_f,[t,x],[dt,dx],[1,0] | 
inhomo=1);
-- nd_weyl_gr :0sec(0.000526sec)
-- weyl_minipoly :0sec(0.0002439sec)
-- generic_bfct_and_gr :0sec(0.001016sec)
generic bfct : [[1,1],[s,1]]
S0 : 0
B_{S0} length : 1
-- fctr(BF) + base :0sec + gc : 0.008sec(0.00691sec)
-- integration_ideal_internal :0sec(0.0003109sec)
[[27*x^3*dx^2+54*x^2*dx+6*x+1],[[[[dt,-dt^2-3*dt-3]],1]]]
\end{verbatim}
}
\end{example}

\begin{theorem} \label{thm:multi-int}
We consider the following multiple integral,  
\begin{equation} \label{eq:multi-int}
F(x_{m+1}, \ldots, x_n) = \int_{a_1}^{b_1} \cdots \int_{a_m}^{b_m} 
  f(x_1, \ldots, x_n) dx_1 \cdots dx_m ~~~ (m \leq n).
\end{equation}
Let $I$ be a holonomic left $D$-ideal annihilating the integrand $f(x_1, \ldots, x_n)$.  
There exists an algorithm to compute inhomogeneous differential equations 
for the multiple integral $F(x_{m+1}, \ldots, x_n)$ from the holonomic ideal $I$.
The algorithm is described below.
\end{theorem}

For simplicity, we will explain the algorithm in the case of $m=2$.
We set 
$$F(x_{3}, \ldots, x_n) = \int_{a_1}^{b_1} \int_{a_2}^{b_2} 
  f(x_1, \ldots, x_n) dx_1 dx_2 ~~~ (2 \leq n), $$
and will compute an inhomogeneous differential equation of $F$. 

Let $I$ be a holonomic left $D$-ideal annihilating the integrand 
$f(x_1, \ldots, x_n)$.
We compute the integration ideal $J$ of $I$ with respect to $x_1, x_2$, 
i.e. 
$$ J = (I + \pd{1} D + \pd{2} D) \cap D' ~~~ 
(D' = K \langle x_3, \ldots, x_n, \pd{3}, \ldots, \pd{n} \rangle).$$
We take an element $P \in J$. 
There exist $P_0 \in I$ and $P_1, P_2 \in D$ 
such that $P = P_0 + \pd{1} P_1 + \pd{2} P_2 \in D'.$ 
We apply the operator $P$ to the integral $F$, and obtain
\begin{equation}
P \cdot F = \int_{a_2}^{b_2} (P_1 \cdot f|_{x_1 = b_1} - P_1 \cdot f|_{x_1 = a_1})dx_2
    + \int_{a_1}^{b_1}(P_2 \cdot f|_{x_2 = b_2} - P_2 \cdot f|_{x_2 = a_2}) dx_1.
\label{mulint:1}
\end{equation}
Let $F_1, F_2$ be the first term and the second term of the right hand side 
and let $f_1, f_2$ be the integrand of $F_1, F_2$. 

To obtain a holonomic ideal annihilating the integral $F_1$, 
we must compute a holonomic ideal $I_1$ annihilating the integrand $f_1$.  
When the integrand $f_1$ is the power of polynomial, 
we can use Oaku's algorithm to obtain the holonomic ideal $I_1$ (\cite{OAKU}).
In general case, 
we can compute the holonomic ideal $I_1$ from $I$ by the following method. 

The ideal quotient $I : P_1$ is holonomic and annihilates the function $P_1 \cdot f$.
To obtain a holonomic ideal $J_1$ annihilating $P_1 \cdot f|_{x_1 = b_1}$,
we compute the restriction ideal of $I : P_1$ with respect to $x_1=b_1$.
Applying the same procedure for $x_1 = a_1$ instead of $x_1 = b_1$, 
we obtain a holonomic ideal $J_2$ annihilating 
$P_1 \cdot f|_{x_1 = a_1}$.
Since $J_1 \cap J_2$ is holonomic and annihilates 
$f_1 (= P_1 \cdot f|_{x_1 = b_1} - P_1 \cdot f|_{x_1 = a_1})$, 
we obtain $J_1 \cap J_2$ as $I_1$. 

We compute the integration ideal $K_1$ of $I_1$ with respect to $x_2$, i.e.  
$$K_1 = (I_1 + \pd{2}D_1) \cap D' ~~~ 
(D_1 = K \langle x_2, x_3, \ldots, x_n, \pd{2},\pd{3}, \ldots, \pd{n} \rangle).$$
We take an element $P^{(1)} \in K_1$.
There exist $P_0^{(1)} \in I_1$ and $P_2^{(1)} \in D_1$ 
such that $P^{(1)} = P_0^{(1)} + \pd{2} P_2^{(1)}$.
We apply $P^{(1)}$ to the integral $F_1$, and obtain
\begin{equation}
P^{(1)} \cdot F_1 = P_2^{(1)} \cdot f_1|_{x_2 = b_2} - P_2^{(1)} \cdot f_1|_{x_2 = a_2}. \label{mulint:2}
\end{equation}
Applying the same procedure for $I_2$ instead of $I_1$, 
we can compute the annihilating ideal $I_2$ of the integrand $f_2$  
and the integration ideal $K_2$ of $I_2$ with respect to $x_1$. 

By (\ref{mulint:1}) and (\ref{mulint:2}), 
we obtain 
$$ P^{(1)} \cdot P \cdot F = P^{(1)} \cdot F_1 + P^{(1)} \cdot F_2, $$
and can compute the first term of the right hand side.
To compute the second term $P^{(1)} \cdot F_2$,
we compute $K_2 : P^{(1)}$ and take an element $P^{(2)}$ in this ideal.
Since $P^{(2)} P^{(1)} \in K_2$, we can compute $P^{(2)} P^{(1)} \cdot F_2$.
Finally, we can obtain an inhomogeneous differential equation 
$$ P^{(2)} P^{(1)} P \cdot F = P^{(2)} P^{(1)} \cdot F_1 + P^{(2)} P^{(1)} \cdot F_2.$$

\begin{remark} \normalfont
Let 
$$\ell_1 \cdot F = g_1, \cdots, \ell_p \cdot F = g_p \quad 
(\ell_i \in D',~ g_i \text{ is a holonomic function })$$
be a system of inhomogeneous differential equations.
When $\langle \ell_1, \ldots, \ell_p \rangle$ generates the left holonomic ideal in $D'$, 
we call the system {\it inhomogeneous holonomic}.
When $m = 1$, the output of the algorithm in Theorem \ref{thm:multi-int} 
is inhomogeneous holonomic.
Although the algorithm outputs a lot of inhomogeneous differential equations 
when $P$ runs over the ideal $J$, 
it is an open question whether the output of the algorithm is inhomogeneous holonomic when $m > 1$.
However, since the Oaku-Shiraki-Takayama algorithm gives holonomic output
(see \cite{OST}, section \ref{sec:OST}), 
we can obtain inhomogeneous holonomic differential equations by the following algorithm.
\end{remark}

\begin{algorithm}
~
\normalfont
\begin{itemize}
\item[Input:] Generators of a holonomic left ideal annihilating
$f(x_1, \ldots, x_n)$.
\item[Output:]  Generators of an inhomogeneous holonomic system
for (\ref{eq:multi-int}).
\end{itemize}
\begin{enumerate}
\item Apply the algorithm in Theorem \ref{thm:multi-int}.
\item Apply the Oaku-Shiraki-Takayama algorithm if the system obtained
in the step 1 is not inhomogeneous holonomic.
\item Merge the outputs of the step 1 and the step 2.
\end{enumerate}
\end{algorithm}

\section{Comparison of our algorithm with other algorithms}
\subsection{The Almkvist-Zeilberger algorithm}
The Almkvist-Zeilberger algorithm (AZ algorithm, \cite{AZ}, \cite{Tefera}, \cite{APZ}) 
is very fast, but works for hyperexponential functions.
Our algorithm works for holonomic functions.  
The AZ algorithm is based on the method of undetermined coefficients 
and Gosper's algorithm, 
and our algorithm is based on the Gr\"{o}bner basis method in $D$. 

\subsection{The Chyzak algorithm}
The Chyzak algorithm (\cite{C1}, \cite{C2}, \cite{CS}) is based on the method of 
undetermined coefficients and the Gr\"{o}bner basis method in the Ore algebra.
By using the Ore algebra, the Chyzak algorithm can compute various summations and integrals
like summations of holonomic sequences, integrals of holonomic functions and its $q$-analogues. 
For the ring of differential operators with rational function 
coefficients $K(x)\langle \partial \rangle$, 
the Chyzak algorithm is a generalization of the AZ algorithm and works for holonomic functions.
The algorithm is often faster than our algorithm.   
But, when the algorithm returns higher order differential equations or the number of variables are many, 
our algorithm is sometimes faster.
Here, we show only one example.
We present these examples at {\tt http://www.math.kobe-u.ac.jp/OpenXM/Math/i-hg/nk\_restriction\_ex.html}

\begin{example} \normalfont
[$F(x,y)=\int^b_a \frac{1}{xt+y+t^{10}} dt$]

We set 
$$ F(x,y) = \int^b_a \frac{1}{xt+y+t^{10}} dt.$$
We will compute differential equations for the integral $F(x,y)$.
The following output is computed by our algorithm.
It takes about 1.3 seconds.
{\small
\begin{verbatim}
[2345] load("nk_restriction.rr");
[2545] F=x*t+y+t^10$             
[2546] Ann=ann(F)$ /* annihilating ideal of F^s */
0.052sec(0.0485sec)
[2547] Id=map(subst, Ann, s, -1)$ /* substitute s=-1 in Ann */
0sec(4.411e-05sec)
[1569] nk_restriction.integration_ideal(Id,[t,x,y],[dt,dx,dy],[1,0,0]
|inhomo=1);
-- nd_weyl_gr :0.012sec + gc : 0.008001sec(0.02009sec)
-- weyl_minipoly :0sec(0.001189sec)
-- generic_bfct_and_gr :0.016sec + gc : 0.008001sec(0.02358sec)
generic bfct : [[1,1],[s,1],[s-9,1]]
S0 : 9
B_{S0} length : 10
-- fctr(BF) + base :0.044sec + gc : 0.024sec(0.0674sec)
-- integration_ideal_internal :0.8321sec + gc : 0.236sec(1.071sec)
[[9*x*dx+10*y*dy+9,-10*dx^9-x*dy^9,-9*dx^10+y*dy^10+9*dy^9],
[[[[dt,-t]],1],[[[dt,-dy^8]],1],[[[dt,-t*dy^9]],1]]]
0.9081sec + gc : 0.28sec(1.19sec)
\end{verbatim}
}

The following output is computed by the Chyzak algorithm 
(package {\tt Mgfun} \cite{Mgfun}) on Maple12.
It takes about 50 seconds.
{\small
\begin{verbatim}
with(Mgfun):
f:=1/(x*t+y+t^10):
ts:=time():
creative_telescoping(f,[x::diff,y::diff], t::diff):
time()-ts;
                                   49.583
\end{verbatim}
}

These computational experiments are executed on a Linux machine with Intel Xeon5450 (3.00GHz) and 
32 GB memory.

\end{example}

\subsection{The Oaku-Shiraki-Takayama algorithm} \label{sec:OST}
Although our algorithm gives inhomogeneous differential equations for definite integrals, 
the Oaku-Shiraki-Takayama algorithm (OST algorithm, \cite{OST}) is for
computing homogeneous differential equations annihilating definite
integrals by using the Heaviside function and the integration algorithm.
Since outputs are different, they are different methods.
However, in most examples, outputs of our algorithm can be easily transformed to homogeneous systems.
Thus, it will be worth making comparison between our method and the OST method.

Let $u(t,x)$ be a smooth function defined on an open neighborhood of 
$[a, b] \times U $ where $U$ is an open set of ${\bf R}^{n-1}$.
The Heaviside function $Y(t)$ defined by 
$Y(t)=0 \: (t<0), Y(t)=1 \: (t \geq 0)$.
Then we can regard the integral of $u(t,x)$ over $[a, b]$ as
that of $Y(t-a)Y(b-t)u(t,x)$ over $(-\infty, \infty)$,
and the following holds.
$$ \int_a^b u(t,x) dt = \int_{-\infty}^{\infty} Y(t-a)Y(b-t)u(t,x)dt $$
Thus we can apply Algorithm \ref{alg:int-ideal} 
to obtain homogeneous differential equations.
The paper \cite{OST} proposes the two methods 
\begin{enumerate}
 \item[(a)] Method of using properties of the Heaviside function
 \item[(b)] Method of using tensor product in $D$-module
\end{enumerate}
to obtain differential equations annihilating the integrand of the right hand side. 
In the former case, the computation finishes without a heavy part
because the procedure is only multiplication of polynomials. 
However, it is not known whether the output is holonomic.
In the latter case, when an input is holonomic,
an output is also holonomic.
However, the computation is often heavy.
We call the former OST algorithm (a) and 
the latter OST algorithm (b) in this paper.
See \cite[Chap 5]{OST} for details.

Let us show a relation of the outputs of OST algorithm and our algorithm.
We consider $v(x) = \int_0^\infty e^{(-t^3+t)x} dt$.
OST algorithm (a) or (b) return the following ideal
\begin{align*}
    \langle & -27x^3 \pd{x}^3-54x^2 \pd{x}^2+(4x^3+3x)\pd{x}+4x^2-3,\\
    & 27x^2 \pd{x}^4+135x \pd{x}^3+(-4x^2+105) \pd{x}^2-16x \pd{x}-8 \rangle .
\end{align*}
On the other hand, Algorithm \ref{alg:int-inhomo} returns
the following ideal generated by $P$ and its inhomogeneous part $Q$:
\begin{align*}
\langle P \rangle &= \langle -27x^2 \pd{x}^2-27x \pd{x}+4x^2+3 \rangle, \\ 
                Q &= \pd{t}(-9tx\pd{x}+(-6t^2+4)x+3t).
\end{align*}
This output yields 
\begin{equation} \label{eq:OST-inhomo}
 P \cdot v(x) 
 = \left [ (-9tx\pd{x}+(-6t^2+4)x+3t) \cdot e^{(-t^3+t)x} \right ]_{t=0}^{t=\infty}
 = -4x.
\end{equation}
Since the annihilating ideal of $-4x$ is
$\langle x \pd{x}-1, \pd{x}^2 \rangle$,
operators $(x \pd{x}-1)  P$ and $\pd{x}^2 P$ annihilate $v(x)$.
Although results of these algorithms are not coincide in general,
these operators coincide outputs of OST algorithm (a) and (b) in this case.
However, it seems that it is difficult to compute the right hand side 
of (\ref{eq:OST-inhomo}) from the output of OST algorithm.
Moreover, in our algorithm we have only to do substitution process
to compute for the integrals which has same integrand and 
another integration domain
because our algorithm does not depend on the integration domain.

Table \ref{table:comp_ost} shows the computing time of each part of
Algorithm \ref{alg:int-inhomo} and OST algorithm (a), (b).
The entries with parentheses for inputs $\bar{v}_k$ 
mean that results for $v_k$ were reused.
For a comparison we show the computing time of Algorithm \ref{alg:int-ideal}.
The experiments were done on a Linux machine
with Intel Xeon X5570 (2.93GHz) and 48 GB memory.
\begin{table}[htbp]
\begin{center}
\begin{tabular}{|c||c|c|c||c||c|c|c|||c|} \hline
  & \multicolumn{3}{c||}{Alg \ref{alg:int-inhomo}} & OST (a) &
 \multicolumn{3}{c|||}{OST (b)} & Alg \ref{alg:int-ideal} \\ \hline
Input & \makebox[3em]{Alg \ref{alg:int-inhomo}} & \makebox[3em]{{\rm
 Ann}} & \makebox[3em]{Total} & \makebox[3em]{Total} &
 \makebox[3em]{(b)} & \makebox[3em]{Alg \ref{alg:int-ideal}} & \makebox[3em]{Total} & \makebox[3em]{Total} \\ \hline
$v_1$       & 0.0042 & 0.0014 & 0.0056 & 0.0062  & 0.11  & 0.012   & 0.12    & 0.0039 \\ \hline
$v_2$       & 0.15   & 0.019  & 0.17   & 0.25    & 5.10  & 0.16    & 5.26    & 0.075  \\ \hline
$v_3$       & 19.91  & 0.45   & 20.36  & 96.14   & 24.54 & 95.24   & 119.8   & 13.58  \\ \hline
$v_4$       & 26724  & 28.33  & 26752  &$>$ 1 day& 1726  &$>$ 1 day& ---     & 24003  \\ \hline \hline

$\bar{v}_1$ &(0.0042)& 0.0015 & 0.0057 & 0.0071   & 0.47  & 0.0050   & 0.48  & n/a \\ \hline
$\bar{v}_2$ &(0.15)  & 0.027  & 0.18   & 1.56     & 18230 & 1.19     & 18231 & n/a \\ \hline
$\bar{v}_3$ &(19.91) & 1.62   & 21.53  & 3769     & 848   & 2802     & 3650  & n/a \\ \hline
$\bar{v}_4$ &(26724) & 294    & 27018  &$>$ 1 day & 16231 &$>$ 1 day & ---   & n/a \\ \hline
\end{tabular}
 $$v_k(x) = \int_0^{\infty} u_k(t,x) dt, \:
 \bar{v}_k(x) = \int_0^{1} u_k(t,x) dt \:\:
 \text{where} \:\:
 \displaystyle u_k(t,x) = \exp \left (-tx \prod_{i=1}^k (t^2-i^2) \right) $$
\end{center}
\caption{The comparison of the computing time (seconds)} \label{table:comp_ost}
\end{table}

From the viewpoint of the computational efficiency,
the computation time of Algorithm \ref{alg:int-inhomo} 
increases more than that of Algorithm \ref{alg:int-ideal} for
computing inhomogeneous parts.
That of OST algorithm increases because the input data
of the integration algorithm becomes bigger differential operators by procedure (a) or (b).
It seems that Algorithm \ref{alg:int-inhomo} is faster than OST algorithm,
since the computation of inhomogeneous parts can be done by
multiplication and summation of differential operators.
However, to obtain homogeneous equation corresponding to OST algorithm
output, we must compute annihilating ideals of inhomogeneous parts.

\section*{Acknowledgement}
We would like to thank Prof. Takayama for fruitful discussions 
and encouragements.

\end{document}